\documentclass{article}
\usepackage{amssymb,amsmath,theorem,euscript}

\newcounter{sec}

\def\sm{}

\newcounter{punct}[sec]

\def\punct{\refstepcounter{punct}{\arabic{sec}.\arabic{punct}.  }}

\def\COUNTERS{\addtocounter{sec}{1}
              \setcounter{punct}{0}
          \setcounter{equation}{0}
                  }

\begin{document}

 \def\ov{\overline}
\def\wt{\widetilde}
 \newcommand{\rk}{\mathop {\mathrm {rk}}\nolimits}
\newcommand{\Aut}{\mathop {\mathrm {Aut}}\nolimits}
\newcommand{\Out}{\mathop {\mathrm {Out}}\nolimits}
\renewcommand{\Re}{\mathop {\mathrm {Re}}\nolimits}
\renewcommand{\Im}{\mathop {\mathrm {Im}}\nolimits}
\def\Br{\mathrm {Br}}

\def\SL{\mathrm {SL}}
\def\SU{\mathrm {SU}}
\def\GL{\mathrm {GL}}
\def\U{\mathrm U}
\def\OO{\mathrm O}
 \def\Sp{\mathrm {Sp}}
 \def\SO{\mathrm {SO}}
\def\SOS{\mathrm {SO}^*}
 \def\Diff{\mathrm{Diff}}
 \def\Vect{\mathfrak{Vect}}
\def\PGL{\mathrm {PGL}}
\def\PU{\mathrm {PU}}
\def\PSL{\mathrm {PSL}}
\def\Symm{\mathrm{Symm}}
\def\Herm{\mathrm{Herm}}
\def\Pos{\mathrm{Pos}}
\def\Lat{\mathrm{Lat}}
\def\Tri{\mathrm{Triang}}

\def\End{\mathrm{End}}
\def\Mor{\mathrm{Mor}}
\def\Aut{\mathrm{Aut}}
 \def\PB{\mathrm{PB}}
 \def\cA{\mathcal A}
\def\cB{\mathcal B}
\def\cC{\mathcal C}
\def\cD{\mathcal D}
\def\cE{\mathcal E}
\def\cF{\mathcal F}
\def\cG{\mathcal G}
\def\cH{\mathcal H}
\def\cJ{\mathcal J}
\def\cI{\mathcal I}
\def\cK{\mathcal K}
 \def\cL{\mathcal L}
\def\cM{\mathcal M}
\def\cN{\mathcal N}
 \def\cO{\mathcal O}
\def\cP{\mathcal P}
\def\cQ{\mathcal Q}
\def\cR{\mathcal R}
\def\cS{\mathcal S}
\def\cT{\mathcal T}
\def\cU{\mathcal U}
\def\cV{\mathcal V}
 \def\cW{\mathcal W}
\def\cX{\mathcal X}
 \def\cY{\mathcal Y}
 \def\cZ{\mathcal Z}
\def\0{{\ov 0}}
 \def\1{{\ov 1}}
 \def\frA{\mathfrak A}
 \def\frB{\mathfrak B}
\def\frC{\mathfrak C}
\def\frD{\mathfrak D}
\def\frE{\mathfrak E}
\def\frF{\mathfrak F}
\def\frG{\mathfrak G}
\def\frH{\mathfrak H}
\def\frI{\mathfrak I}
 \def\frJ{\mathfrak J}
 \def\frK{\mathfrak K}
 \def\frL{\mathfrak L}
\def\frM{\mathfrak M}
 \def\frN{\mathfrak N} \def\frO{\mathfrak O} \def\frP{\mathfrak P} \def\frQ{\mathfrak Q} \def\frR{\mathfrak R}
 \def\frS{\mathfrak S} \def\frT{\mathfrak T} \def\frU{\mathfrak U} \def\frV{\mathfrak V} \def\frW{\mathfrak W}
 \def\frX{\mathfrak X} \def\frY{\mathfrak Y} \def\frZ{\mathfrak Z} \def\fra{\mathfrak a} \def\frb{\mathfrak b}
 \def\frc{\mathfrak c} \def\frd{\mathfrak d} \def\fre{\mathfrak e} \def\frf{\mathfrak f} \def\frg{\mathfrak g}
 \def\frh{\mathfrak h} \def\fri{\mathfrak i} \def\frj{\mathfrak j} \def\frk{\mathfrak k} \def\frl{\mathfrak l}
 \def\frm{\mathfrak m} \def\frn{\mathfrak n} \def\fro{\mathfrak o} \def\frp{\mathfrak p} \def\frq{\mathfrak q}
 \def\frr{\mathfrak r} \def\frs{\mathfrak s} \def\frt{\mathfrak t} \def\fru{\mathfrak u} \def\frv{\mathfrak v}
 \def\frw{\mathfrak w} \def\frx{\mathfrak x} \def\fry{\mathfrak y} \def\frz{\mathfrak z} \def\frsp{\mathfrak{sp}}
 \def\bfa{\mathbf a} \def\bfb{\mathbf b} \def\bfc{\mathbf c} \def\bfd{\mathbf d} \def\bfe{\mathbf e} \def\bff{\mathbf f}
 \def\bfg{\mathbf g} \def\bfh{\mathbf h} \def\bfi{\mathbf i} \def\bfj{\mathbf j} \def\bfk{\mathbf k} \def\bfl{\mathbf l}
 \def\bfm{\mathbf m} \def\bfn{\mathbf n} \def\bfo{\mathbf o} \def\bfp{\mathbf p} \def\bfq{\mathbf q} \def\bfr{\mathbf r}
 \def\bfs{\mathbf s} \def\bft{\mathbf t} \def\bfu{\mathbf u} \def\bfv{\mathbf v} \def\bfw{\mathbf w} \def\bfx{\mathbf x}
 \def\bfy{\mathbf y} \def\bfz{\mathbf z} \def\bfA{\mathbf A} \def\bfB{\mathbf B} \def\bfC{\mathbf C} \def\bfD{\mathbf D}
 \def\bfE{\mathbf E} \def\bfF{\mathbf F} \def\bfG{\mathbf G} \def\bfH{\mathbf H} \def\bfI{\mathbf I} \def\bfJ{\mathbf J}
 \def\bfK{\mathbf K} \def\bfL{\mathbf L} \def\bfM{\mathbf M} \def\bfN{\mathbf N} \def\bfO{\mathbf O} \def\bfP{\mathbf P}
 \def\bfQ{\mathbf Q} \def\bfR{\mathbf R} \def\bfS{\mathbf S} \def\bfT{\mathbf T} \def\bfU{\mathbf U} \def\bfV{\mathbf V}
 \def\bfW{\mathbf W} \def\bfX{\mathbf X} \def\bfY{\mathbf Y} \def\bfZ{\mathbf Z} \def\bfw{\mathbf w}
 \def\R {{\mathbb R }} \def\C {{\mathbb C }} \def\Z{{\mathbb Z}} \def\H{{\mathbb H}} \def\K{{\mathbb K}}
 \def\N{{\mathbb N}} \def\Q{{\mathbb Q}} \def\A{{\mathbb A}} \def\T{\mathbb T} \def\P{\mathbb P} \def\G{\mathbb G}
 \def\bbA{\mathbb A} \def\bbB{\mathbb B} \def\bbD{\mathbb D} \def\bbE{\mathbb E} \def\bbF{\mathbb F} \def\bbG{\mathbb G}
 \def\bbI{\mathbb I} \def\bbJ{\mathbb J} \def\bbL{\mathbb L} \def\bbM{\mathbb M} \def\bbN{\mathbb N} \def\bbO{\mathbb O}
 \def\bbP{\mathbb P} \def\bbQ{\mathbb Q} \def\bbS{\mathbb S} \def\bbT{\mathbb T} \def\bbU{\mathbb U} \def\bbV{\mathbb V}
 \def\bbW{\mathbb W} \def\bbX{\mathbb X} \def\bbY{\mathbb Y} \def\kappa{\varkappa} \def\epsilon{\varepsilon}
 \def\phi{\varphi} \def\le{\leqslant} \def\ge{\geqslant}

\def\UU{\bbU}
\def\Mat{\mathrm{Mat}}
\def\tto{\rightrightarrows}

\def\Gr{\mathrm{Gr}}

\def\graph{\mathrm{graph}}

\def\O{\mathrm{O}}

\def\la{\langle}
\def\ra{\rangle}

\def\B{\mathrm B}
\def\Int{\mathrm {Int}}

\begin{center}

{\bf\Large Matrix beta-integrals: an overview}

\vspace{10pt}

\sc 
 Yurii A. Neretin%
 \footnote{These paper is based on notes of my talk in XXXIII Workshop 'Geometry
 and Physics' Bialowieza, 2014. Supported by the grant FWF, Project P25142}

\end{center}

{\small First examples of matrix beta-integrals were discovered on 1930-50s by
Siegel and Hua, in  60s Gindikin obtained multi-parametric series
of such integrals. We discuss beta-integrals related to symmetric spaces,
their interpolation with respect to the dimension of a ground field, and adelic analogs;
also we discuss beta-integrals related to flag spaces.}

\vspace{22pt}

{\sc
\noindent
1. Introduction. The Euler and Selberg integrals
\\
2. The Hua integrals
\\
3. Beta-functions of symmetric spaces
\\
4. Zeta-functions of spaces of lattices
\\
5. Non-radial interpolation of matrix beta-function
\\
6. Beta-functions of flag spaces
}

\section{Introduction. The Euler and Selberg integrals}

\COUNTERS

{\bf \punct Euler beta-function.}
Recall the standard formulas for the Euler beta-function:
\begin{align}
&\int_0^1 x^{\alpha-1}(1-x)^{\beta-1}dx=
\frac{\Gamma(\alpha)\Gamma(\beta)}{\Gamma(\alpha+\beta)}
\quad &&\text{(Euler)}
\label{eq:beta-1}
\\
&\int_\R \frac{dx}{(1+ix)^\mu(1-ix)^\nu} =
\frac{2^{2-\mu-\nu}\pi\Gamma(\mu+\nu-1)}{\Gamma(\mu)\Gamma(\nu)}
\quad &&\text{(Cauchy)}
\label{eq:beta-2}
\\
&\int_0^\infty \frac{x^{\alpha-1}}{(1+x)^\sigma}=
\frac{\Gamma(\alpha)\Gamma(\sigma-\alpha)}{\Gamma(\sigma)}
&&
\label{eq:beta-3}
\\
&\int_0^\pi (\sin t)^\mu e^{i\nu t}dt=
\frac\pi{2^\mu} \frac{\Gamma(1+\mu)}
{\Gamma(1+\frac{\mu+\nu}2)\Gamma(1+\frac{\mu-\nu}2)}
e^{i\pi\nu/2}
\quad &&\text{(Lobachevsky)}
\label{eq:beta-4}
\end{align}
The integral (\ref{eq:beta-3}) is obtained from  (\ref{eq:beta-1})
by the substitution $x=t/(1+t)$. Replacing the segment $[0,1]$
in (\ref{eq:beta-1}) by the circle $|x|=1$, after simple manipulations we get
(\ref{eq:beta-4}). Considering the stereographic projection of the circle to the line,
we come to (\ref{eq:beta-2}).

\sm

{\bf\punct Beta-integrals.} 'Beta-integral' is an informal term for
integrals of the type
\begin{equation}
\int \Bigl( \text{Product})=
\text{Product of Gamma-functions}
\label{eq:beta-integrals}.
\end{equation}

There is large family of such identities (see, e.g., \cite{Ask}, \cite{AAR}).
 First, we  present
two nice examples.
The De Branges \cite{dB} -- Wilson integral (1972, 1980)
is given by
$$
\frac 1{2\pi}\int\limits_\R
\Bigl|\frac {\prod_{j=1}^4\Gamma(a_j+ix) }
{\Gamma(2ix)}  \Bigr|^2\,\,dx=
\frac{\prod_{1\le k<l\le 4}\Gamma(a_k+a_l)}{\Gamma(a_1+a_2+a_3+a_4)}
.
$$

Recall that the integrand is a weight function for the Wilson orthogonal polynomials,
which occupy the highest level of the Askey hierarchy \cite{KS} of hypergeometric
orthogonal polynomials. 

The second example is the Selberg integral, \cite{Sel}, 1944,
\begin{multline} 
\int_0^1\dots\int_0^1
\prod_{j=1}^n t_j^{\alpha-1}(1-t_j)^{\beta-1}
\prod_{1\le k<l\le n}|t_k-t_l|^{2\gamma}\,
dt_1\dots dt_n=
\\
=\prod_{j=1}^{n}
\frac{\Gamma(\alpha+(j-1)\gamma)\,\Gamma(\beta+(j-1)\gamma)\,
\Gamma\bigl(1+j\gamma\bigr)}
{\Gamma\bigl(\alpha+\beta+(n+j-2)\gamma\bigr)\,\Gamma(1+\gamma)}
.
\label{eq:selberg}
\end{multline}
As the Euler beta-integral, the Selberg integral has several versions, for instance
\begin{multline}
\int_0^\infty\dots\int_0^\infty
\prod_{j=1}^n x_j^{\alpha-1}(1+x_j)^{-\alpha-\beta-2\gamma(n-1)}
\prod_{1\le k<l\le n}|x_k-x_l|^{2\gamma}\,
dx_1\dots dx_n
=
\\
=\prod_{j=1}^{n}
\frac{\Gamma(\alpha+(j-1)\gamma)\,\Gamma(\beta+(j-1)\gamma)\,
\Gamma\bigl(1+j\gamma\bigr)}
{\Gamma\bigl(\alpha+\beta+(n+j-2)\gamma\bigr)\,\Gamma(1+\gamma)}
\label{eq:selberg2}
\end{multline}
\noindent
\begin{multline}
\frac 1{(2\pi)^n}
\int\limits_{-\infty}^\infty\dots \int\limits_{-\infty}^\infty
\prod_{k=1}^n(1-ix_k)^{-\alpha}(1+ix_k)^{-\beta}
\prod_{1\le k<l\le n}|x_k-x_l|^{2\gamma}\,dx_1\dots dx_n
=\\=
2^{-(\alpha+\beta)n+\gamma n(n-1)+n}
\prod_{j=1}^n
\frac{\Gamma\bigl(\alpha+\beta  -(n+j-2)\gamma-1\bigr)\Gamma(1+j\gamma)}
{\Gamma\bigl(\alpha-(j-1)\gamma)\bigr)\Gamma\bigl(\beta-(j-1)\gamma)\bigr)
 \Gamma(1+\gamma) }
\label{eq:selberg1}
\end{multline}

There exists a large family of beta-integrals (\ref{eq:beta-integrals}), including one-dimensional
integrals  (see an old 
overview of Askey \cite{Ask}), multi-dimensional
integrals, $q$-analogs, elliptic analogs; some occasional 
collection
of references is \cite{AAR}, \cite{KS}, \cite{Koo}, \cite{RS}, \cite{Ner-wil}, \cite{Spi},
\cite{FW}.

The topic of these notes is
analogs of integrals 
(\ref{eq:beta-1})--(\ref{eq:beta-4}), (\ref{eq:selberg})--(\ref{eq:selberg1}).

\sm

{\bf\punct Notation.}

$\bullet$ $\K$ denote $\R$, $\C$, or quaternions $\H$,
$
\frd:=\dim \K.
$

\sm

$\bullet$ $[X]_p$ is the left upper corner of a matrix $X$ of size $p\times p$;

$\bullet$ $[X]_{pq}$ is the left upper corner of a matrix $X$ of size $p\times q$;

$\bullet$ $X^*$, $X^t$ are the adjoint matrix and the transposed matrix;

$\bullet$  $X>0$ means that a matrix $X$ is self-adjoint and {\it strictly} positive definite,
$X>Y$ means that $X-Y>0$;

$\bullet$ $\|X\|$ denotes a {\it norm of a matrix}, precisely the norm of the corresponding linear
 operator in the standard
Euclidean space.
, $\|X\|=\|X^* X\|^{1/2}= \|X X^*\|^{1/2}$; for a self-adjoint matrix
norm is $\max |\lambda_j|$ over all eigenvalues.

\sm

Spaces of matrices:

\sm
 
$\bullet$ $\Mat_{p,q}(\K)$ is the space of all matrices of size $p\times q$ over $\K$;

$\bullet$ $\Herm_n(\K)$ is the space of all Hermitian matrices ($X=X^*$) of size $n$;

$\bullet$ $\Symm_n(\K)$ is the space of all symmetric matrices ($X=X^t$)  of size $n$.


\sm

The Lebesgue measure on such spaces is normalized in the most simple way.
For instance,  for $\Symm_n(\R)$ we set
$$
dX:=\prod\limits_{1\le k\le l \le n} dx_{kl};
$$
for $\Mat_{p,q}(\C)$, we write
$$
dZ:= \prod\limits_{1\le k\le p,\,1\le l \le q} d\Re z_{kl}\, d\Im z_{kl}.
$$

\section{The Hua integrals}

\COUNTERS

{\bf\punct The Hua integrals.} The famous book  \cite{Hua}
'{\it Harmonic analysis of functions of several complex variables in classical domains}'
by Hua, 1958,
 contains calculations of a family of matrix  integrals. We present two
 examples.
 
Consider the space $\B_{m,n}$ of complex $m\times n$ matrices $Z$
with $\|Z\|<1$. The following identity holds
\begin{equation}
\int\limits_{ZZ^*<1} \det(1-ZZ^*)^\lambda\,dZ=
\frac{\prod_{j=1}^n \Gamma(\lambda+j) \prod_{j=1}^m \Gamma(\lambda+j)}
{\prod_{j=1}^{n+m} \Gamma(\lambda+j)}\pi^{nm}
.
\label{eq:hua-1}
\end{equation}

Next, consider the space $\Symm_n(\R)$ of all real symmetric matrices of size $n$.
The following identity holds
\begin{equation}
\int\limits_{\Symm_n(\R)}
\frac{dT}{\det(1+T^2)^\alpha}
=
\pi^{\frac{n(n+1)}4}
\frac{\Gamma(\alpha-n/2)}{\Gamma(\alpha)}
\prod_{j=1}^{n-1} \frac{\Gamma\bigl(2\alpha-(n+j)/2\bigr)}
{\Gamma(2\alpha-j)} 
.
\label{eq:hua-2}
\end{equation}

{\bf\punct Comments: spaces and integrands.} 
We can consider the following 10 series of {\it matrix spaces}

\sm

$\bullet$ $p\times q$ matrices over $\R$;

$\bullet$ symmetric $n\times n$ matrices ($X=X^t$) over $\R$;

$\bullet$ skew-symmetric $n\times n$ matrices ($X=-X^t$) over $\R$;

$\bullet$ $p\times q$ matrices over $\C$;

$\bullet$ symmetric $n\times n$ matrices  over $\C$;

$\bullet$ skew-symmetric $n\times n$ matrices  over $\C$;

$\bullet$ Hermitian $n\times n$ matrices ($X=X^*$) over $\C$;

$\bullet$ $p\times q$ matrices over $\H$;

$\bullet$ Hermitian $n\times n$ matrices ($X=X^*$) over $\H$;

$\bullet$ anti-Hermitian $n\times n$ matrices ($X=-X^*$) over $\H$.

\sm

For any space of this list, we consider a '{\it matrix ball}' $XX^*<1$.

\sm

For all 'matrix spaces' and all 'matrix ball',  integrals%
\footnote{Recall a definition of a {\it determinant $\det (X)=\det_\H (X)$ of a quaternionic matrix} $X$.
Such matrix determines a transformation $\H^n\to\H^n$ and therefore an $\R$-linear transformation
$X_\R:\R^4\to\R^4$. We set $\det_\H (X):=\sqrt[4]{\det(X_\R)}$. In particular,
$\det(X)$ is real non-negative. 
If entries of $X$ are complex, then the quaternionic
determinant coincides with $|\det_\C X|$.
}
\begin{align}
&\int \det(1+XX^*)^{-\alpha} dX ;
\label{eq:hua-3}
\\
& \int_{XX^*<1} \det(1-XX^*)^{\gamma} dX 
 \label{eq:hua-4}
\end{align}
are long products of gamma-functions as (\ref{eq:hua-1})--(\ref{eq:hua-2}).
Actually, Hua evaluated 1/3 of these  20 integrals. 
Apparently, there is no text, where all these integrals are evaluated 
(and a reason, which does not excuse this,  is explained in the next subsection).

The domain of integration $\B_{m,n}\subset \C^{nm}$  in
(\ref{eq:hua-1}), i.e.,
the  matrix ball $\|Z\|<1$, is a well-known object in differential
geometry,  representation theory, and complex analysis, since it is 
an Hermitian 
 symmetric space%
 \footnote{Spaces $\B_{p,q}$ also are known as {\it Cartan domains} of type I.},
 $$
 \B_{p,q}=\U_{p,q}\bigr/(\U_p\times \U_q)
$$
The pseudounitary group $\U_{p,q}$ acts on this domain by
linear-fractional transformations:
\begin{equation}
\begin{pmatrix}
a&b\\c&d
\end{pmatrix}:Z\mapsto U:= (a+Zc)^{-1}(b+Zd)
.
\label{eq:lf}
\end{equation}

The remaining 9 series of 'matrix balls' $XX^*<1$ also are Riemannian symmetric spaces%
\footnote{below a '{\it symmetric space}' means a semisimple (reductive) 
 symmetric space.}.
Up to a minor inaccuracy, all Riemannian noncompact symmetric spaces
admit  'matrix ball' models. The group of isometries consists of 
certain linear-fractional transformations (see tables of symmetric spaces 
in \cite{Ner-gauss}, Addendum D).

Meaning of the integrand $\det(1-ZZ^*)^\alpha$ is less obvious%
\footnote{Hua Loo Keng evaluated volumes of Cartan domains and some compact symmetric spaces
and observed that calculations survive in a wider generality.}.
However, any mathematician what had deal with the unit circle $|z|<1$ could observe
that the expression $(1-z\ov z)^\alpha$ quite often appears in formulas.
The same holds for  $\det(1-ZZ^*)^\alpha$ in the case of the matrix balls.
We only point out a nice behavior
of the expression under linear-fractional transformation (\ref{eq:lf}):
$$
\det(1-UU^*)^\alpha= \det(1-ZZ^*)^\alpha |\det(a+zc)|^{-2\alpha}
.$$
Thus integrals (\ref{eq:hua-4}) are integrals of some reasonable expressions
over  non-compact symmetric spaces. 

Integrals (\ref{eq:hua-3})
are integrals over compact symmetric spaces written in   coordinates.
For instance, in (\ref{eq:hua-2}) we integrate over the space $\Symm_n(\R)$.
But $\Symm_n(\R)$ is a chart on the real Lagrangian Grassmannian
(recall that if an operator $T:\R^n\to\R^n$ is symmetric, then
its graph is a Lagrangian subspace in $\R^n\oplus\R^n$, see, e.g.,
\cite{Ner-gauss}, Sect.3.1). The Lagrangian Grassmannian is 
a homogeneous (symmetric) space $\U_n/\OO_n$, see, e.g.,
\cite{Ner-gauss}, Sect. 3.3. All other 'matrix spaces' defined above
are open dense charts on certain compact Riemannian symmetric spaces.
Up to a minor inaccuracy, all compact symmetric spaces
admit such charts (see tables of symmetric spaces 
in \cite{Ner-gauss}, Addendum D).

\sm

{\bf\punct Integration over eigenvalues.}
Consider the space $\Herm_n(\K)$ of all Hermitian matrices%
\footnote{$\Herm_n(\R)$ is $\Symm_n(\R)$.}
 over $\K=\R$, $\C$, or $\H$;
equip this space with the standard Lebesgue measure.
To a matrix $X\in\Herm_n(\K)$, we assign the collection of its eigenvalues
\begin{equation}
\Lambda:\,\lambda_1\ge\lambda_2\ge\dots\ge \lambda_n.
\label{eq:lambda}
\end{equation}
Thus we get a map $X\mapsto\Lambda$ from $\Herm_n(\K)$ to the wedge
(\ref{eq:lambda}). The distribution of eigenvalues is given by the formula
$$
C_n(\K)
\prod\limits_{1\le k<l\le n}|\lambda_k-\lambda_l|^{\frd} d\lambda_1\dots d\lambda_n
,$$
where $C_n(\K)$ is a certain (explicit) constant, $\frd=\dim K$. This can be reformulated as follows.
Let 
$F$  be a function on $\Herm_n(\K)$ invariant with respect 
 to the unitary group $\U(n,\K)$%
\footnote{$\U(n,\R)$ is the orthogonal group $\OO(n)$, $\U(n,\C)$ is the usual unitary
group $\U(n)$, $\U(n,\H)$ is the compact symplectic group $\Sp(2n)$.},
$$
F(uXu^{-1})= F(X),\qquad u\in\U(n,\K)
.$$
Such $F$ is a function of eigenvalues,
$$
 F(X)=f(\lambda_1,\dots,\lambda_n).
$$
Then  the following integration formula holds
\begin{multline}
\int\limits_{\Herm_n(\K)} F(X)\,dX= \\= C_n(\K)
\int\limits_{\lambda_1\ge\lambda_2\ge\dots\ge \lambda_n}
f(\lambda_1,\dots,\lambda_n)
\prod_{1\le k<l\le n}|\lambda_k-\lambda_l|^{\frd} d\lambda_1\dots d\lambda_n
.
\label{eq:weyl}
\end{multline}
The formula is  a relative of the Weyl integration formula,
see derivations of several formulas of this kind in \cite{Hua}.
   
In the Hua integral (\ref{eq:hua-2}), the integrand is
$$
\det(1+T^2)^{-\alpha}=\prod\limits_{j=1}^n (1+\lambda_j^2)^{-\alpha}=
\prod\limits_{j=1}^n (1+i\lambda_j)^{-\alpha} (1-i\lambda_j)^{-\alpha}
.$$
Applying the integration formula (\ref{eq:weyl}) we reduce the Hua 
integral (\ref{eq:hua-2}) to a special case of the Selberg integral
(\ref{eq:selberg1}). Moreover, we get also an explicit evaluation of 
a more general integral
$$
\int
\det (1+iT)^\alpha \det(1-iT)^\beta\,dT
.$$

Next, consider the space of all complex matrices of size $m\times n$, where
$m\le n$.
To each  matrix we assign a collection of its singular values%
\footnote{{\it Singular values} of a matrix $Z$ are eigenvalues of $\sqrt{ZZ^*}$.}
$$\mu_1\ge\mu_2\ge\dots\ge \mu_m\ge 0.$$
The distribution of singular values is given by 
$$
\prod\limits_{1\le k\le n} \mu_k^{2(n-m)+1}
\prod\limits_{1\le k<l\le m} (\mu_k^2-\mu_l^2)^2 
\prod\limits_{1\le k\le n} d\mu_k.
$$
The integrand in the Hua integral (\ref{eq:hua-1}) is
$\prod (1-\mu_k^2)^2$. After the substitution $x_k=\mu_k^2$, this integral also is reduced
to the Selberg integral (\ref{eq:selberg1}).

All  20 integrals (\ref{eq:hua-3})--(\ref{eq:hua-4})
are reduced to the Selberg integrals in a similar
 way%
\footnote{In all these cases the parameter  $\gamma$ in the Selberg integrals
is 1/2, 1, 2. For some exceptional symmetric spaces distributions of invariants
give $\gamma=4$.}.

\sm

{\bf\punct An application of Hua calculations: projective systems of measures.}
Let us return to  integral (\ref{eq:hua-2}). Represent a matrix $T$ as a block
matrix of size
$(n-1)+1$,
$$
T=\begin{pmatrix} S&p\\p^t&q
\end{pmatrix}
.$$
Consider a function $f$ on $\Symm_n(\R)$ depending only on $S=[T]_{n-1}$.
 Then the following identity
holds
\begin{multline}
\int\limits_{\Symm_n(\R)} f(S)\det\biggl(1+\begin{pmatrix} S&p\\p^t&q
\end{pmatrix}^2\biggr)^{-\alpha}\,dS\, dp\,dq
=\\=
2^{\frac{n-1}2}\pi^{\frac n2}
\frac{\Gamma(2\alpha+\frac{n+1}2)\Gamma(\alpha-\frac12)}
{\Gamma(\alpha)\Gamma(2\alpha-1)}
\int\limits_{\Symm_{n-1}(\R)} f(S)\det(1+S^2)^{1/2-\alpha}dS
.
\label{eq:projectivity}
\end{multline}
This formula can be extracted from the original Hua calculation
(the formula (\ref{eq:projectivity}) also implies (\ref{eq:hua-2})). 

Now fix $\alpha>-1/2$ and consider a measure $\nu_{\alpha,n}$ on $\Symm_n$
given by 
$$\nu_{\alpha,n}=
s_{\alpha,n} \det(1+T^2)^{-\alpha-(n+1)/2}dT
,$$
where the normalizing constant $s_{\alpha,n}$ is chosen to make the
total  measure $=1$. Consider the chain of projections
$$
\dots \longleftarrow \Symm_{n-1}(\R)
 \longleftarrow \Symm_{n}(\R) \longleftarrow \dots
,$$
where each map sends a matrix 
$X\in\Symm_n(\R)$ to its left upper corner $[X]_{n-1}(\R)$. According
(\ref{eq:projectivity}), this map sends the measure $s_{\alpha,n}$
to the measure $s_{\alpha,n-1}$. By the Kolmogorov
consistency theorem (see, e.g., \cite{Shi}, \S2.9)  
 there is a measure $\nu_\alpha$
on the space $\Symm_{\infty}(\R)$ of infinite symmetric matrices
whose image under each map $X\mapsto [X]_n$ is $\nu_{\alpha,n}$.

Next, consider the group of finitary%
\footnote{We say that a matrix $g$ is {\it finitary}, if $g-1$
has finite number of nonzero matrix elements} orthogonal $(\infty+\infty)$
block matrices  having the  structure $\begin{pmatrix}a&b\\-b&a\end{pmatrix}$.
 This group is isomorphic to the group $\U_\infty$
  of finitary unitary matrices. It acts on  $\Symm_{\infty}(\R)$
  by linear-fractional transformations (\ref{eq:lf}), point out that  this formula makes sense. It is easy to show that the measure $\nu_\alpha$ is quasiinvariant with respect to such transformations, and there arises a  problem of decomposition
  of the space $L^2$. We also can regard our limit space as the
  inverse limit of the chain of Lagrangian Grassmannians,
  $$
  \dots \longleftarrow \U_{n-1}/\OO_{n-1}
 \longleftarrow \U_{n}/\OO_{n} \longleftarrow \dots
  $$
  
Such construction exists for any series of compact symmetric spaces
and leads to an interesting harmonic analysis on the limit objects,
see \cite{Pic}, \cite{Ner-hua}, \cite{Ols}, \cite{BO}.

\sm

{\bf\punct Remarks.}
 1) The construction of inverse limits does not admit an extension to non-compact symmetric
spaces (i.e., to matrix balls). Of course, the chain of projections of
{\it sets}
$$
\dots \longleftarrow \B_{p,q}
\longleftarrow \B_{p+1,q+1}
\longleftarrow \B_{p+2,q+2}
\longleftarrow
\dots
$$
is well defined. We can consider normalized probabilistic measures
$$
s'_{\alpha,p,q,k} \det (1-ZZ^*)^{\alpha-2k}
$$
on $\B_{p+k,q+k}$. However, for sufficiently large $k$ the integral
$$
\int_{\B_{p+k,q+k}}\det (1-ZZ^*)^{\alpha-2k}\,dZ
$$
is divergent.

\sm

2) Projective limits exist for $p$-adic Grassmannians, see \cite{Ner-p}.

\section{Beta-functions of symmetric spaces}

\COUNTERS

{\bf \punct The Gindikin beta-function of symmetric cones.} 
Consider the space $\Pos_n(\K)$  of positive definite $n\times n$
matrices over $\K$. The cone $\Pos_n(\K)$ is
a model of the symmetric space $\GL_n(\K)/\U_n(\K)$, the group
$\GL_n(\K)$ acts on $\Pos_n(\K)$ by  transformations
$$
g:\, X\mapsto g^* Xg.
$$

Gindikin \cite{Gin}, 1965, considered a matrix $\Gamma$-function given by
\begin{multline}
\mathbf\Gamma[\mathbf{s}]:=\int_{\Pos_n(\K)}e^{-\mathrm{tr} \, X}
\prod_{j=1}^n \det [X]_j^{s_j-s_{j+1}}\cdot
\det X^{\frd n/2- \frd/2+1}\,dX
=\\=(2\pi)^{n(n-1)\frd/4}
\prod_{k=1}^n \Gamma\Bigl(s_k- (k-1)\frac \frd 2
\Bigr)
.
\label{eq:gin-gamma}
\end{multline}
Here $s_j\in\C$, $s_{n+1}:=0$; $[X]_p$ denotes upper left corners
of size $p$ of a matrix $X$.
 The expressions $s_j-s_{j+1}$
are written by aesthetic reasons, we can write 
$$\prod\nolimits_{j=1}^n \det [X]_j^{\lambda_j}$$
with arbitrary $\lambda_j$.
The factor $\det X^{\frd n/2- \frd/2+1}$ can be included to the latter
product, but it is the density of the $\GL_n(\K)$-invariant measure
on $\Pos_n(\K)$ and it is reasonable to split it from the product.

To evaluate the integral, Gindikin considers%
\footnote{See also, \cite{FK}.} the substitution $X=S^*S$,
where $S$ is an upper triangular matrix with positive elements on the diagonal.
After this the integral splits into a product of one-dimensional integrals.

Also the following imitation of the beta-function:
take place
\begin{multline}
\mathbf{B}[\mathbf{s},\mathbf{t}]:=
\int\limits_{0<X<1}
\prod_{j=1}^n\Bigl(\det [X]_j^{s_j-s_{j+1}}
\cdot \det [1-X]_j^{t_j-t_{j+1}}\Bigr)
\times\\ \times
\det X^{\frd n/2- \frd/2+1} \det (1-X)^{\frd n/2- \frd/2+1}\, dX
=
\frac{\mathbf\Gamma[\mathbf{s}]\,\mathbf\Gamma[\mathbf{t}]}
{\mathbf\Gamma[\mathbf{s}+\mathbf{t}]}
.
\label{eq:gin-beta}
\end{multline}
A proof in \cite{Gin} is an one-to-one imitation of the standard evaluation
of the Euler beta-integral.

These integrals extend some results of 1920-30s
(Whishart, Ingham, Siegel, see \cite{Sie}).

\sm

{\bf\punct Beta functions of Riemannian non-compact symmetric spaces.}
The domain of integration $0<X<1$ in 
(\ref{eq:gin-beta}) is itself  the symmetric space 
$\GL_n(\K)/\U_n(\K)$. Indeed,
the matrix ball $ZZ^*<1$ in the space of Hermitian matrices
is a model of the symmetric space 
$\GL_n(\K)/\U_n(\K)$). The inequality $ZZ^*<1$ is equivalent
to $-1<Z<1$, and we substitute $Z=-1+2X$.

Analogs of integrals (\ref{eq:gin-beta}) for 7 remaining
series of  Riemannian non-compact symmetric spaces were obtained
in \cite{Ner-beta}%
\footnote{For the case of tubes $\SO_0(n,2)/\SO(n)\times \SO(2)$, which
is slightly exceptional, see \cite{Ner-tube}.}.
  We give two well-representative examples.

In the first example we consider a symmetric space, which can be realized
 as a matrix wedge.
Let $W_n$ be the  domain  ({\it Siegel upper-half plane}) 
of $n\times n$ complex symmetric matrices $Z$ with
 $\Re Z>0$.
This is a model of a symmetric space $\Sp_{2n}(\R)/\U_n$.
We write $Z=T+iS$, where $T$, $S$ are real symmetric matrices. Then
\begin{multline}
\int\limits_{T=T^t>0,\, S=S^t}
\prod_{j=1}^n
\frac{ \det[T]_j^{\lambda_j-\lambda_{j+1}}    }
     { \det[1+T+iS]_j^{\sigma_j-\sigma_{j+1}} 
     \det[1+T-iS]_j^{\tau_j-\tau_{j+1}} }
     \times\\\times
    \det T^{-(n+1)} dT\,dS=  \\
= \prod_{k=1}^n
\frac{2^{2-\sigma_k-\tau_k+n-k}  \pi^{k}
 \Gamma(\lambda_k-(n+k)/2)\Gamma(\sigma_k+\tau_k-\lambda_k - (n-k)/2)  }
 {  \Gamma(\sigma_k-(n-k)/2) \Gamma(\tau_k-(n-k)/2)}
  \label{eq:wedge}
\end{multline}
(we set $\lambda_{j+1}=\sigma_{j+1}=\tau_{j+1}=0$).

There are also noncompact symmetric spaces, which do not admit realizations
as convex matrix cones and convex matrix wedges. As an example,
 we consider the space $\OO_{p,q}/\OO_{p}\times \OO_{q}$.
Let $q\ge p$. We realize this space (for details, see \cite{Ner-beta}, Sect.3)
 as the space of real 
block matrices of size $(q-p)+p$
having the form
$$R=
\begin{pmatrix}
1&0\\2L&K
\end{pmatrix}
$$
and satisfying the dissipativity condition 
$$
R+R^t>0
$$
We represent $K$ as $K=M+N$, where $M$ is is symmetric and $K$ is skew-symmetric.
Then the dissipativity condition $R+R^t>0$ reduces to the form
$$
\begin{pmatrix}
1&L^t\\L&M
\end{pmatrix}>0
$$
or equivalently $M-LL^t>0$. We have the following integrals in 
coordinates $L$, $M$, $K$:
\begin{multline}
\int\limits_{\footnotesize\begin{matrix}
M=M^t>0,N=-N^t\\
M-LL^t>0
\end{matrix}}
\prod_{j=1}^p
\frac{ \det[M-LL^t]_j^{\lambda_j-\lambda_{j+1}}  }
  {\det [1+M+N]_j^{\sigma_j-\sigma_{j+1}} }
  \times\\\times
\det(M-LL^t)^{-(p+q)/2} dM\,dN\,dL =\\
=\prod_{k=1}^p
\pi^{k-(q-p)/2-1}
\frac{\Gamma(\lambda_k-(q+k)/2+1)\Gamma(\sigma_k-\lambda_k-(p-k)/2)}
     {\Gamma(\sigma_k-p+k)}
     .
     \label{eq:section-wedge}
     \end{multline}

{\bf \punct Remarks.}
1) Integrals (\ref{eq:wedge})-({\ref{eq:section-wedge})
 were written to obtain Plancherel measure for Berezin representations
of classical groups, see \cite{Ner-beta}, \cite{Ner-Plan}.

\sm

2) I do not know perfect counterparts of the integrals
(\ref{eq:wedge})-({\ref{eq:section-wedge}) for  compact symmetric spaces.
Some beta-integrals over classical groups $\SO(n)$, $\U(n)$, $\Sp(n)$
were considered in \cite{Ner-hua}, extensions to over compact symmetric spaces
 are more-or-less automatic. However, they depend on a smaller number of parameters.
 
 \sm
 
3) {\it On analogs of the $\Gamma$-function.} To be definite, consider the space 
$\Mat_{n,n}(\C)$. Consider a distribution
$$
\phi(Z)=
\prod\nolimits_{j=1}^n |\det [Z]_j|^{\lambda_j} \det [Z]_j^{p_j},
$$
where $p_j\in\Z$, $\lambda_j\in \C$. This expression is homogeneous
in the following sense: for an upper triangular matrix $A$ and a lower triangular matrix
$B$,
$$
\phi(BZA)=\prod |a_{jj}b_{jj}|^{\sum_{k\le j} \lambda_j} 
(a_{jj}b_{jj})^{\sum_{k\le j} p_j} \phi(Z).
$$
The Fourier transform $\widehat \phi$ of $\phi$ must be homogeneous. For $\lambda_j$
in a general position this remark allows to write $\widehat \phi$ up to a constant factor.
This factor (it is a product of Gamma-functions and sines) can be regarded
as a matrix analog of Gamma-function. See Stein \cite{Ste}, 1967, Sato, Shintani
\cite{SS}, 1974. I do not know an exhausting text on this topic.

\section{Zeta-functions of spaces of lattices}

\COUNTERS

Noncompact symmetric spaces have $p$-adic counterparts, namely Bruhat--Tits buildings
(see, e.g., \cite{Ner-gauss}, Chapter 10).
Since this topic is not inside common knowledge, we will discuss an adelic
variant of matrix beta-integrals. 

\sm

{\bf\punct Space of lattices.}
A {\it lattice} in $\Q^n$ is a subgroup isomorphic to $\Z^n$. 
Denote by $\Lat_n$   the space of lattices in $\Q^n$.
The group $\GL_{n}(\Q)$ acts on the space  $\Lat_n$,
the stabilizer of the standard lattice $\Z^n$ is $\GL_n(\Z)$. Thus
$\Lat_n$ is a homogeneous space
$$
\Lat_n\simeq \GL_n(\Q)/\GL_n(\Z).
$$
 
{\bf\punct Analog of beta-integrals.}
We consider two coordinate flags
\begin{align*}
0\subset \Z\subset \Z^2\subset\dots\subset \Z^n;
\\
0\subset \Q\subset \Q^2\subset\dots\subset \Q^n
.
\end{align*}
Consider intersections of a lattice $S$ with these flags,
i.e.,
$$
S\cap \Z^k\,\,\subset\,\, S\cap  \Q^k\,\,\subset\,\, \R^k
.
$$
For a lattice $S\subset \R^k$ we denote by
 $\upsilon_k(S)$  the volume of the quotient $\R^k/S$.
 The following identity  holds \cite{Ner-tits}:
\begin{multline}
\sum\limits_{S\in\mathrm{Lat}_n(\Q)}\,\,
\prod\limits_{j=1}^n
\upsilon_k(S\cap\Q^k)^{-\beta_k+\beta_{k+1}}
\,\,
\upsilon_k(S\cap\Z^k)^{-\alpha_k+\alpha_{k+1}}
=\\=
\prod\limits_{j=1}^n
\frac{\zeta(-(\beta_j+j-1))\,\,\zeta(\alpha_j+\beta_j-n+j)}
{\zeta(\alpha_j-n+j)}
\label{eq:lat}
,\end{multline}
where
$\zeta$ is the Riemannian $\zeta$-function,
$$
\zeta(s)=\sum_{k=0}^\infty \frac 1{n^s}=
\prod_{\text {prime $p$}} \Bigl(1-\frac 1{p^s} \Bigr)^{-1}
.
$$

{\bf \punct On Berezin kernels.} It seems that holomorphic discrete
series  representations 
of semisimple Lie groups have no $p$-adic analogs. However,
in \cite{Ner-tits} there were obtained analogs of the Berezin kernels and of
the Berezin--Wallach set.  Let us explain this on our minimal language.
We define a Berezin kernel on $\Lat_n$ by
$$
K_\alpha(S,T):=\frac{\bigl(\upsilon_n(R)\,\upsilon_n(S)\bigr)^{\alpha/2}}
{\bigl(\upsilon_n(R\cap S)\bigr)^{\alpha}}
.
$$
This kernel is positive definite if and only if
$$
\alpha=0,1,\dots, n-1,\text{or $\alpha>n-1$} 
.
 $$
Positive definiteness of the kernel means that there exists a Hilbert space
$H_\alpha$ and a total system of vectors $\delta_S\in H_\alpha$,
 where $S$ ranges in $\Lat_n$, such that
$$
\la \delta_S,\delta_T\ra_{H_\alpha}=K_\alpha(S,T)
.
$$
The group $\GL_n(\Q)$ acts in the spaces $H_\alpha$. Further picture
is parallel to the theory of Berezin kernels over $\R$
(see \cite{Ner-Plan}).
 Formula (\ref{eq:lat}) allows to obtain
the Plancherel formula for this representation.

\sm

{\bf\punct Remarks.}
1) An analog of $\Gamma$-function is the  Tamagawa zeta-function \cite{Tam}, see also \cite{Mac}. It  is the  sum
$$
\sum \prod\nolimits_{k=1}^n \upsilon_k(S\cap\Z^k)^{-\alpha_k+\alpha_{k+1}}
$$
 over  sublattices is $\Z^n$. It can be obtained from
(\ref{eq:lat}) by a degeneration.

\sm

b) Certainly, analogs of (\ref{eq:lat}) for symplectic and orthogonal groups
must exist. As far as I know they are not yet obtained.

\section{Non-radial interpolation  of matrix beta-integrals}

\COUNTERS

{\bf\punct Rayleigh tables.}
Again, $\K=\R$, $\C$, or quaternions $\H$, $\frd=\dim \K$.
Consider Hermitian matrices of order $n$ over $\K$.

Consider eigenvalues of $[X]_p$ for each $p$,
$$
\lambda_{p1}\le\lambda_{p2}\le\dots\le \lambda_{pp}.
$$
We get a  table $\cL$  

\begin{equation}
\begin{array} {ccccccccccc}
 & & & &   & \lambda_{11}  &  & & & &  \\
 & & & &  \lambda_{21}  & & \lambda_{22}  & & & &  \\
 & & &\lambda_{31} &   &\lambda_{32} &  &\lambda_{33} & & &  \\
 & &\cdot &\cdot & \cdot  &\cdot & \cdot &\cdot &\cdot & &  \\
 &\cdot &\cdot &\cdot & \cdot  &\cdot & \cdot &\cdot &\cdot &\cdot &  \\
\lambda_{n1} & &\lambda_{n2} & &\lambda_{n3}
      &\dots & \lambda_{n(n-2)} & &\lambda_{n(n-1)}
        &  & \lambda_{nn}
\end{array}
\label{eq:table}
\end{equation}
with the {\it Rayleigh interlacing condition}%
\footnote{This statement also is called the Rayleigh--Courant-Fisher theorem.}
$$
\dots\le \lambda_{(j+1)k}\le \lambda_{jk}\le \lambda_{(j+1)(k+1)}\le \dots
$$
This means that the numbers $\lambda_{kl}$ increase in 'north-east' and
 'south-east' directions.

Denote by $\cR_n$ the space of all {\it Rayleigh tables} (\ref{eq:table}).

Point out that for $\K=\R$ the number of variables $\lambda_{kl}$ coincides
with $\dim \Herm_n(\R)$ (but generally there are $2^{n(n-1)/2}$ matrices
$X$ with a given $\cL$).

Now consider the image of the Lebesgue measure on $\Herm_n(\K)$
under the map $\Herm_n(\K)\to \cR_n$. In other words, consider the joint distribution
of eigenvalues of all $[X]_p$. It is given by the formula
\begin{multline}
d\rho_\frd(\cL) = C_n(\frd)
\frac
{\prod\limits_{2\le j\le n} \quad
\prod\limits_{1\le \alpha\le j-1,\,\,\,
                               1\le p\le j}
|\lambda_{(j-1)\alpha}-\lambda_{jp}|^{\frd/2-1} }
{\prod\limits_{2\le j\le n-1}\quad
\prod\limits_{1\le \alpha<\beta\le j}
  (\lambda_{j\beta}-\lambda_{j\alpha})^{\frd-2}}
\times \\ \times
\prod\limits_{1\le p<q\le n}
(\lambda_{nq} -\lambda_{np})
 \prod\limits_{1\le j\le n}\quad \prod_{1\le\alpha\le j}
  d\lambda_{j\alpha}
  ,
  \label{eq:measure}
\end{multline}
where
$$
C_n(\frd)=\frac{\pi^{n(n-1)\frd/4}}
{\Gamma^{n(n-1)/2} (\frd/2)}
.
$$
Notice that for $\K=\C$ we get a total cancellation in the expression
(\ref{eq:measure}).
History of this formula is not quite clear. It seems that ideologically
it is contained  in  book \cite{GN} by Gelfand, Naimark (see evaluation
of spherical functions of $\GL(n,\C)$). The measure
(\ref{eq:measure}) is used in integral representation of Jack polynomials
in paper \cite{OO} by Olshanski and Okounkov. A formal proof is contained in \cite{Ner-Ray}, see also \cite{FR} and \cite{Kaz2}.

\sm

{\bf\punct Interpolation.}
Now we can assume that $\frd$ is an arbitrary complex number and interpolate
matrix beta-integrals 
\begin{multline*}
\int_{\Herm_n(\K)}\prod_{k=1}^{n-1} (1+i[X]_k)^{-\sigma_k+\sigma_{k+1}-\frd/2}
(1-i[X]_k)^{-\tau_k+\tau_{k+1}-\frd/2}\times\\ \times
 \det(1+iX)^{-\sigma_n}
\det(1-iX)^{-\tau_n} dX=\frac{\prod\Gamma(\dots)}{\prod\Gamma(\dots)}
\end{multline*}
 with respect to $\frd=\dim \K$:
\begin{multline*}
\int_{\cR_n}
\prod\limits_{j=1}^{n-1}\prod\limits_{\alpha=1}^j
(1+i\lambda_{j\alpha})^{-\sigma_j+\sigma_{j+1}-\frd/2}
(1-i\lambda_{j\alpha})^{-\tau_j+\tau_{j+1}-\frd/2}
\times\\ \times \prod\limits_{p=1}^n
(1+i\lambda_{np})^{-\sigma_n}
(1-i\lambda_{np})^{-\tau_n}
 d\rho_\frd(\Lambda)
=\\=
\pi^{n(n-1)\frd/4+n}
\cdot
\prod\limits_{j=1}^n
\frac{\Gamma(\sigma_j+\tau_j-1-(j-1) \frd/2)}
    {\Gamma(\sigma_j)\Gamma(\tau_j)}
    .
\end{multline*}
Here integration is taken over the space of all Rayleigh tables
and the measure
$d\rho_\frd(\Lambda)$ is given by (\ref{eq:measure}).

However, the proof \cite{Ner-Ray} of the latter formula remains to be valid 
for a wider family of integrals,
\begin{multline*}
\int
\prod\limits_{j=1}^{n-1}\prod\limits_{\alpha=1}^j
(1+i\lambda_{j\alpha})^{-\sigma_j+\sigma_{j+1}-\theta_{j\alpha}}
(1-i\lambda_{j\alpha})^{-\tau_j+\tau_{j+1}-\theta_{j\alpha}}
\times\\ \times \prod\limits_{p=1}^n
(1+i\lambda_{np})^{-\sigma_n}
(1-i\lambda_{np})^{-\tau_n}
\times\\ \times
\prod\limits_{j=1}^{n-1}
\frac{
  \prod\limits_{1\le\alpha\le j,\, 1\le p\le j+1}
   |\lambda_{j\alpha}-\lambda_{(j+1)p}|^{\theta_{j\alpha}-1}
   }
   {
  \prod\limits_{1\le\alpha<\beta \le j}
  (\lambda_{j\beta}-\lambda_{j\alpha})^{\theta_{j\alpha}+\theta_{j\beta}-2}
  }
\prod\limits_{1\le p< q\le n} (\lambda_{nq}-\lambda_{np})\, d\Lambda
=\\=
\pi^{n} 2^{2n-\sum\limits_{j=1}^n(\sigma_j+\tau_j)}
\prod\limits_{1\le\alpha\le j\le n-1} \Gamma(\theta_{j\alpha})
\cdot
\prod\limits_{j=1}^n
\frac{\Gamma(\sigma_j+\tau_j-1-\sum_{\alpha=1}^{j-1} \theta_{(j-1)\alpha})}
    {\Gamma(\sigma_j)\Gamma(\tau_j)}
    .
\end{multline*}
Now the parameter $\frd$ is replaced by $(n-1)n/2$ parameters $\theta_{j\alpha}$

\sm

{\bf\punct Remarks.}
The Gindikin beta-integrals admit an interpolation in the same spirit
\cite{Ner-Ray}. For beta-integrals (\ref{eq:wedge})--(\ref{eq:section-wedge})
over wedges and more general domains an interpolation is unknown. 

\section{Beta-integrals over flag spaces}

\COUNTERS

{\bf\punct Beta-integrals.}
Now we consider upper-triangular matrices $Z=\{z_{ij}\}$ over $\K$, 
$$
z_{ii}=1,\quad z_{ij}=0\quad \text{for $i>j$}.
$$
Denote the space of all upper-triangular matrices by $\Tri_n(\K)$. 
Recall that the space of upper-triangular matrices is a chart on a flag space.

Let $[Z]_{pq}$ be left upper corners of $Z$ of size $p\times q$,
denote
$$
s_{pq}(Z):=\det ([Z]_{pq}[Z]_{pq}^*).
$$
The following identity \cite{Ner-flags} holds
$$
\int_{\Tri_n(\K)}
 \prod_{1\le p<q\le n} s_{pq}(Z)^{-\lambda_{pq}}\,
dZ=
\pi^{n(n-1)/4} \prod_{1\le p<q\le n} \frac{\Gamma(\nu_{pq}-\frd/2)}{\Gamma(\nu_{pq})}
,
$$
where the integration is taken over the space of upper-triangular matrices, and
$$
\nu_{pq}:=-\frac12(q-p-1)\frd+\sum_{k,m:\, p\le k < q,\, q\le m\le n} \lambda_{mk}
.
$$

{\bf\punct Projectivity.}  
 Consider the map $Z\mapsto [Z]_{n-1}$
 from $\Tri_n(\K)$ to $\Tri_{n-1}(\K)$. 
 Consider  a measure%
 \begin{equation*}
 \prod\nolimits_{p=1}^{n-1} s_{pn} (z)^{-\lambda_p}\,dZ^{\{n\}}
 \label{eq:density}
 \end{equation*}
 on $\Tri_n(\K)$. Assume
 $$
 \lambda_{p}+\lambda_{p+1}+\dots +\lambda_{n-1}>\frac 12(n-p)\frd
 $$
 for all $p$.
 Then the pushforward of this measure under the forgetting map
 is 
 \begin{multline*}
 \label{eq:projection}
\pi^{\frac{(n-1)\frd}2}\prod\nolimits_{1\le p\le n-1}
\frac{\Gamma(\lambda_p+\dots +\lambda_n-(n-p)\frd/2)}
{\Gamma(\lambda_p+\dots +\lambda_n-(n-p+1)\frd/2)}
\times\\\times
 \prod\nolimits_{p=1}^{n-2} s_{p(n-1)} ([Z]_{n-1})^{-\lambda_p}\,\,d[Z]_{n-1}
 .
\end{multline*}

\noindent
\tt Math.Dept., University of Vienna,
 \\
 Oskar-Morgenstern-Platz 1, 1090 Wien;
 \\
\& Institute for Theoretical and Experimental Physics (Moscow);
\\
\& Mech.Math.Dept., Moscow State University.
\\
e-mail: neretin(at) mccme.ru
\\
URL:www.mat.univie.ac.at/$\sim$neretin
} 

\end{document}